\documentclass[11pt,reqno]{amsart}

\usepackage[latin1]{inputenc}
\usepackage[left=3.5cm,right=3.5cm,top=3.5cm,bottom=3.5cm]{geometry}
\usepackage{amssymb}
\usepackage{graphicx}
\usepackage{amscd}
\usepackage[hidelinks]{hyperref}
\usepackage{color}
\usepackage{float}
\usepackage{graphics,amsmath,amssymb}
\usepackage{amsthm}
\usepackage{amsfonts}
\usepackage{latexsym}
\usepackage{epsf}
\usepackage{enumerate}
\usepackage{xifthen}
\usepackage{mathrsfs}
\usepackage{dsfont}
\usepackage{makecell}
\usepackage[FIGTOPCAP]{subfigure}
\usepackage{amsmath}
\allowdisplaybreaks[4]
\usepackage{listings}
\usepackage{etoolbox}
\usepackage{fancyhdr}

 \newtheoremstyle{mytheorem}% name % cf. thmtest.tex of AMSLaTeX
 {3pt}%      Space above
 {3pt}%      Space below
 {\slshape}% Body font
 {}%         Indent amount (empty = no indent,
 {\bfseries}% Thm head font
 {.}%        Punctuation after thm head
 { }%        Space after thm head: " " = normal interword space;
 {}%         Thm head spec (can be left empty, meaning `normal')

\numberwithin{equation}{section}

\theoremstyle{mytheorem}
\newtheorem{theorem}{Theorem}[section]

\newtheorem{lemma}[theorem]{Lemma}

\theoremstyle{definition}

\newcommand{\Keywords}[1]{\ifthenelse{\isempty{#1}}{}{\smallskip \smallskip \noindent \textbf{Keywords}. #1}}
\newcommand{\MSC}[2][2010]{\ifthenelse{\isempty{#2}}{}{\smallskip \smallskip \noindent \textbf{#1MSC}. #2}}
\newcommand{\abstractnote}[1]{\ifthenelse{\isempty{#1}}{}{\smallskip \smallskip \noindent \textsuperscript{\dag}#1}}

\makeatletter
\def\specialsection{\@startsection{section}{1}%
  \z@{\linespacing\@plus\linespacing}{.5\linespacing}%
  {\normalfont}}% NEW
\def\section{\@startsection{section}{1}%
  \z@{.7\linespacing\@plus\linespacing}{.5\linespacing}%
  {\normalfont\scshape}}% NEW
\patchcmd{\@settitle}{\uppercasenonmath\@title}{\Large}{}{}
\patchcmd{\@settitle}{\begin{center}}{\begin{flushleft}}{}{}
\patchcmd{\@settitle}{\end{center}}{\end{flushleft}}{}{}
\patchcmd{\@setauthors}{\MakeUppercase}{\normalsize}{}{}
\patchcmd{\@setauthors}{\centering}{\raggedright}{}{}
\patchcmd{\section}{\scshape}{\large\bfseries}{}{}
\renewcommand{\@secnumfont}{\bfseries}
\patchcmd{\@startsection}{\@afterindenttrue}{\@afterindentfalse}{}{}
\patchcmd{\abstract}{\leftmargin3pc}{\leftmargin1pc}{}{}

\def\maketitle{\par
  \@topnum\z@ % this prevents figures from falling at the top of page 1
  \@setcopyright
  \thispagestyle{empty}% this sets first page specifications
  \ifx\@empty\shortauthors \let\shortauthors\shorttitle
  \else \andify\shortauthors
  \fi
  \@maketitle@hook
  \begingroup
  \@maketitle
  \toks@\@xp{\shortauthors}\@temptokena\@xp{\shorttitle}%
  \toks4{\def\\{ \ignorespaces}}% defend against questionable usage
  \edef\@tempa{%
    \@nx\markboth{\the\toks4
      \@nx\MakeUppercase{\the\toks@}}{\the\@temptokena}}%
  \@tempa
  \endgroup
  \c@footnote\z@
  \@cleartopmattertags
}
\makeatother

\title[Integral right triangle and rhombus pairs]{Integral right triangle and rhombus pairs with a common area and a common perimeter}

\author[S. Chern]{Shane Chern}
\address{School of Mathematical Sciences, Zhejiang University, Hangzhou, 310027, China}
\email{\href{mailto:shanechern@zju.edu.cn}{shanechern@zju.edu.cn}; \href{mailto:chenxiaohang92@gmail.com}{chenxiaohang92@gmail.com}}

\date{}

\begin{document}

{\footnotesize\noindent \textit{Forum Geom.} \textbf{16} (2016), 25--27.}

\bigskip \bigskip

\maketitle

\begin{abstract}
We prove that there are infinitely many integral right triangle and $\theta$-integral rhombus pairs with a common area and a common perimeter by the theory of elliptic curves.


\MSC{Primary 51M25; Secondary 11D25.}
\end{abstract}

\section{Introduction}

We say that a polygon is \textit{integral} (resp. \textit{rational}) if the lengths of its sides are all integers (resp. rational numbers). In a recent paper, Y. Zhang \cite{Zhang2016} proved that there are infinitely many integral right triangle and parallelogram pairs with a common area and a common perimeter. This type of problem originates from a question of B. Sands, which asked for examples of such right triangle and rectangle pair; see the paper of R. K. Guy \cite{Guy1995}. Actually, R. K. Guy gave a negative answer to B. Sands' question, whereas in the same paper showed that there are infinitely many such isosceles triangle and rectangle pairs. Later in 2006, A. Bremner and R. K. Guy \cite{BG2006} replaced isosceles triangle by Heron triangle and proved that such pairs are also infinite.

In this note, we consider such right triangle and rhombus pairs with more restrictions. We say that an integral (resp. rational) rhombus is \textit{$\theta$-integral} (resp. \textit{$\theta$-rational}) if both $\sin\theta$ and $\cos\theta$ are rational numbers. Our result is

\begin{theorem}\label{th:main}
There are infinitely many integral right triangle and $\theta$-integral rhombus pairs with a common area and a common perimeter.
\end{theorem}

\section{Proof of the theorem}

We start from rational right triangles and $\theta$-rational rhombi. Without loss of generality, we may assume that the rational right triangle has sides $(1-u^2,2u,1+u^2)$ with $0<u<1$, and the $\theta$-rational rhombus has side $p$ and intersection angle $\theta$ with $0<\theta \le \pi/2$. Here $u$ and $p$ are both positive rational numbers. Now if the right triangle and rhombus have a common area and a common perimeter, then we have the following Diophantine system
\begin{equation}\label{eq:main01}
\begin{cases}
u(1-u^2)=p^2\sin\theta,\\
1+u=2p.
\end{cases}
\end{equation}
Since both $\sin\theta$ and $\cos\theta$ are rational numbers, we may set
$$\sin\theta=\frac{2v}{1+v^2},$$
where $0<v\le 1$ is a rational number. Note that the case $v=1$, that is $\theta=\pi/2$, was studied by R. K. Guy in \cite{Guy1995}. We thus only need to consider cases $0<v<1$. Eliminating $p$ in \eqref{eq:main01}, we have
\begin{equation}\label{eq:main02}
2u^2v^2-2uv^2+2u^2+uv-2u+v=0.
\end{equation}
One readily notices that if \eqref{eq:main02} has infinitely many rational solutions $(u,v)$ with $0<u,v<1$, then there exist infinitely many pairs of rational right triangle and $\theta$-rational rhombus with a common area and a common perimeter, and thus infinitely many such ($\theta$-)integral pairs by the homogeneity of these sides.

Now by the transformation
$$(x,y)=\left(-\frac{4uv^2+4u+4v-4}{v^2},-\frac{8uv^2-4v^2+8u+8v-8}{v^3}\right),$$
and
$$(u,v)=\left(-\frac{x^3+4x^2+2xy-y^2+4x+4y}{4x^2+y^2+16x+16},\frac{2x+4}{y}\right),$$
we deduce the following elliptic curve
\begin{equation}
E:y^2-3xy-12y=x^3+6x^2+8x.
\end{equation}
Through \textit{Magma}, we compute that $E(\mathbb{Q})$ has rank $1$ (generated by point $P=(0,0)$) and a torsion point of order two $T=(-4,0)$. Note that
$$[4]P=\left(\frac{5920}{4761},\frac{5576768}{328509}\right)$$
leads to a solution
$$(u,v)=\left(\frac{552}{1105},\frac{483}{1264}\right)$$
to \eqref{eq:main02} satisfying $0<u,v<1$. This solution immediately gives a right triangle with side lengths being $1832642$, $2439840$, and $3051458$, and a rhombus with side length being $1830985$ and smaller intersection angle being $\arcsin (1221024/1830985)$. They have a common area $2235676628640$ and a common perimeter $7323940$.

At last, recalling the following result due to H. Poincar\'e and A. Hurwitz (\cite[p. 173]{Poi1901}; see also \cite[Satz 13]{Hur1917}):

\begin{lemma}[Poincar\'e-Hurwitz]\label{le:PH}
Let $E$ be a nonsingular cubic curve in $\mathbb{P}^2$ which is defined over $\mathbb{Q}$. If the set $E(\mathbb{Q})$ is infinite, then every open subset of $\mathbb{P}^2(\mathbb{R})$ which contains one point of $E(\mathbb{Q})$ must contain infinitely many points of $E(\mathbb{Q})$.
\end{lemma}

Now since $[4]P$ gives a suitable solution to \eqref{eq:main02}, from the map $(x,y)\mapsto (u,v)$ and Lemma \ref{le:PH}, we conclude that \eqref{eq:main02} has infinitely many solutions $(u,v)$ satisfying $0<u,v<1$. Thus, there are infinitely many pairs of integral right triangle and $\theta$-integral rhombus with a common area and a common perimeter. This ends the proof of Theorem \ref{th:main}.

\bibliographystyle{amsplain}

\end{document}